\documentclass{article}
\usepackage[leqno]{amsmath}
\usepackage{amsthm,amscd}
\usepackage{amsfonts} 
\usepackage{amssymb} 
\newtheorem{theorem}{Theorem}[section]
\newtheorem{definition}[theorem]{Definition}

\newtheorem{lemma}[theorem]{Lemma}
\newtheorem{proposition}[theorem]{Proposition}

\newtheorem{remark}[theorem]{Remark}

\newcommand{\A}{{\mathcal A}}

\newcommand{\sD}{{\mathsf D}}

\newcommand{\Der}{{\rm Der}}

\newcommand{\pppk}{\frac{\partial}{\partial P_k}}
\newcommand{\JP}{J({\mathbf P})}

\begin{document}
\title{The Hodge filtration and the contact-order filtration of
derivations of Coxeter arrangements}
\author{{\sc Hiroaki Terao}
\footnote{partially supported by the Grant-in-aid for scientific research
(No. 14340018 and 13874005), the Ministry of Education, Sports,
Science and Technology, Japan}\\
{\small \it Tokyo Metropolitan University, Mathematics Department}\\
{\small \it Minami-Ohsawa, Hachioji, Tokyo 192-0397, Japan}
}
\date{}
\maketitle

\begin{abstract}

\noindent
The Hodge filtration of the module of derivations on 
the orbit space of a finite real reflection
group acting on an $\ell$-dimensional Euclidean space was introduced 
and studied by K. Saito \cite{sai2} \cite{sai4}.  
The filtration is equivalent data to the
flat structure or the Frobenius manifold
structure. We will show that the Hodge filtration coincides with the filtration
by the order of contacts to the reflecting hyperplanes.
Moreover, a standard basis for the Hodge filtration is
explicitly given. 
\end{abstract}

\bigskip

\setcounter{section}{0}
\setcounter{equation}{0}

\section{Introduction and main results}

Let 
$V$ be an $\ell$-dimensional Euclidean vector space with 
inner product $I$.
Let $W$ be a
finite irreducible orthogonal reflection group 
(a Coxeter group) acting on $V$.
Let ${\mathcal A}$ be the Coxeter arrangement:
$\A$ is the set of reflecting hyperplanes.
Denote the dual vector space of $V$ by $V^{*}$. 
Then  $V^{*}$ is equipped with the inner product $I^{*} $
which is induced from $I$.
Let $S$  be the symmetric algebra of $V^{*} $ over $\mathbb R$,
which is identified with the algebra of polynomial functions on
$V$.  
The algebra $S$ is naturally graded by $S=\oplus_{q \geq 0} S_q $ where
$S_q$ is the space of homogeneous polynomials of degree $q$.
(Then $S_1=V^*$.)
Let ${\rm Der}_S$ be the $S$-module of ${\mathbb R}$-derivations of $S$.
We say that $\theta \in {\rm Der}_S$ is homogeneous of degree $q$ if
$\theta (S_1) \subseteq S_q$.
Choose for each hyperplane $H \in \cal A$ a linear form $\alpha_H
\in V^* $ such that $H=\ker (\alpha_H).$
Let
\begin{eqnarray}
\label{101} 
\sD^{(m)}({\mathcal A})=\{\theta \in {\rm Der}_S ~|~ \theta (\alpha_H) \in 
S\alpha_{H}^{m} {\rm{~for~ any~}} H \in {\cal A} \}
\end{eqnarray}
for each nonnegative integer $m$.
Elements of $\sD^{(m)}(\cal A)$ are called $m$-{\bf derivations} which
were introduced by G. Ziegler \cite {zie}.
Then one has the {\bf contact-order filtration}
of ${\rm Der}_S$:
$$
{\rm Der}_S = \sD^{(0)}({\mathcal A}) \supset \sD^{(1)}({\mathcal A}) 
\supset \sD^{(2)}({\mathcal A}) 
\supset \cdots .
$$
In \cite{ter}, we proved that each
$\sD^{(m)}(\cal A)$ is a free $S$-module
of rank $\ell$. 

The Coxeter 
group $W$ naturally acts on $V^{*}$, $S$  and
$\Der_{S}$.  
The $W$-invariant subring of $S$
is denoted by  $R$.  
Then it is classically known \cite[V.5.3, Theorem 3]{bou} that there
exist algebraically independent homogeneous polynomials $P_1, \cdots, 
P_{\ell}
\in R$ with $\deg P_{1} \leq\dots\leq \deg P_{\ell}$,
which 
are called {\bf basic invariants},
 such that $R={\mathbb R}[P_1, \cdots, P_{\ell}].$
The {\bf primitive derivation} $D \in {\rm Der}_R$ is characterized
 by
\begin{eqnarray*}
D P_i=\begin{cases}
                                        1 & \text{ for } i=\ell,  \\
                                        0 & \text{ otherwise.}
     \end{cases}
\end{eqnarray*}

Let $\Delta$ be an anti-invariant, in other words, a constant multiple of
$Q := \prod_{H\in \A} \alpha_{H}$.  Then $\Delta^{2} \in R$. 
Let
$$I : \Der_{R} \times \Der_{R} \to \frac{1}{\Delta^2} R$$ 
be the symmetric $R$-bilinear form induced from $I$.
Let
$$\nabla:
{\rm Der}_{R} \times {\rm Der}_{R}  \longrightarrow  
\frac{1}{\Delta^{2}}{\rm Der}_{R}$$
$$({ X}, { Y})  \longmapsto  \bigtriangledown_
{ X}{ Y}$$
be
the {\bf Levi-Civita connection} with respect to $I$. 
Define
\[
T := \{f\in R \mid Df = 0\} = {\mathbb R}[P_{1}, \dots ,P_{\ell-1}].
\]
Then $\nabla_{D} : \Der_{R} \longrightarrow \frac{1}{\Delta^{2}}{\rm Der}_{R}$
is a $T$-linear covariant derivative.
  
    In \cite{ter}, we introduced an $S$-basis  
$\xi^{(m)}_1, \dots, \xi^{(m)}_\ell \in \Der_S$
    for $\sD^{(m)}(\A)$ (Definition~\ref{formula})
    for $m\geq 0$. 
Recall, also from \cite{ter}, the
matrices $B^{(k)} ~~~(k\geq 1)$  with
entries in $T$ whose determinants are
nonzero constants (Lemma \ref{2.1}).
Then

\begin{theorem}
\label{2.4}
For $k \geq 1,$ 

(1)
 $(\nabla_D \xi^{(2k+1)}_1, \dots, \nabla_D 
\xi^{(2k+1)}_\ell)=-(\xi^{(2k-1)}_1, \dots, \xi^{(2k-1)}_\ell)
(B^{(k)})^{-1}B^{(k+1)},$ 
   
   (2)
    $(\nabla^k_D \xi^{(2k-1)}_1, \dots, \nabla^k_D \xi^{(2k-1)}_\ell)
=(-1)^{k-1}(\frac{\partial}{\partial P_1}, \dots, \frac{\partial}
{\partial P_\ell}) B^{(k)}.$
\end{theorem}

Define 
${\mathcal G}_0:=\{ \delta \in \Der_{R} \mid \left[ D, \delta \right]=0 \}$ 
and
${\mathcal G}_k:=\{\delta \in \Der_{R} \mid \nabla^k_D \delta 
\in {\mathcal G}_0 \}$ for $k \geq 1$.
Let $p\geq 0$.
 Put ${\mathcal H}^{(p)} := \bigoplus_{k \geq p} {\mathcal G}_{k}$.
The decreasing filtration 
$$
{\mathcal H}^{(0)} \supset
{\mathcal H}^{(1)} \supset
{\mathcal H}^{(2)} \supset
\dots$$
of ${\mathcal H}^{(0)} = \Der_{R} $ 
is called the {\bf Hodge filtration},   
which was introduced by
K. Saito \cite{sai2} \cite{sai4} to study the primitive integrals. 
The following theorem asserts that the Hodge filtration is
equal to the filtration induced from the contact-order filtration.

\begin{theorem}
\label{1.2} 

(1) ${\mathcal G}_k=\bigoplus^\ell_{j=1} T
\xi^{(2k-1)}_j
   ~~~(k \geq 1),$

(2) ${\mathcal H}^{(p)}=\sD^{(2p-1)}({\mathcal A})
\cap \Der_{R}=\bigoplus^\ell_{j=1} R \xi^{(2p-1)}_j ~~~ (p \geq 1).$
\end{theorem}

\section{Proofs}

We will prove Theorems \ref{2.4} and \ref{1.2} in this section.
  In what follows we use the following notation:\\
   
\noindent
$X_1, \cdots, X_{\ell}$:  a basis for $V^*$\\
$\partial_i := \partial / \partial X_{i} $: 
the partial derivation with respect to 
$X_i~~~(i = 1,\dots, \ell)$\\
$A:=(a^{ij})=(I^{*}(X_{i}, X_{j}))$\\
${\bf X}:=(X_1, \dots, X_\ell)$ \\
${\mathbf P}:=(P_{1}, \cdots, P_{\ell})$\\
$m_{j} := \deg P_{j} -1 ~~(j = 1, 2, \dots , \ell)$\\
$h := m_{\ell} + 1 = \deg P_{\ell} $ : the Coxeter number\\  
$J({\bf g}):=(\frac{\partial g_j}{\partial X_i}):$ the Jacobian matrix for 
${\bf g}=(g_1, \dots, g_\ell)$ \\
$\delta \left[M\right]:=(\delta (m_{ij}))$ for a matrix
$M=(m_{ij})$ and any mapping $\delta$\\
$D^{j} := D\circ D\circ \dots \circ D $ ($j$ times)\\
$\nabla_{D}^{j} := \nabla_{D}\circ \nabla_{D}
\circ \dots \circ \nabla_{D} $ ($j$ times).\\

Let $k \geq 1$. Define
$$
B^{(k)}:=-J({\bf P})^T AJ(D^k[{\bf X}])J(D^{k-1}[{\bf X}])^{-1}J({\bf P})
$$
as in \cite{ter}. 
Then we have

\begin{lemma}
\label{2.1} 
(1) Every entry of $B^{(k)}$  lies in $T$;
$D[B^{(k)}]=0$, 

(2) $\det B^{(k)} \in {\mathbb R}^{*}  $, 

(3) $\deg B^{(k)}_{ij}=m_i +m_j -h,$

(4) $B^{(k+1)}-B^{(k)}=B^{(1)}+(B^{(1)})^T.$ 
\end{lemma}

\begin{proof}
\cite[Lemma 3.2, Lemma 3.4]{ter}. 
\end{proof} 

The  Levi-Civita connection with respect to $I$ 
$$\nabla:
{\rm Der}_{R} \times {\rm Der}_{R}  \longrightarrow  
\frac{1}{\Delta^{2}}{\rm Der}_{R}$$
$$({ X}, { Y})  \longmapsto  \bigtriangledown_
{ X}{ Y}$$
is characterized by the following two properties:\\

(A)  ${ X}(I( Y, Z)) = 
I(\bigtriangledown_{ X}{ Y}, { Z}) + 
I( Y, \bigtriangledown_{ X}{ Z})$  (compatibility),\\
    
(B) $\bigtriangledown_{ X}{ Y} - 
\bigtriangledown_{ Y}{ X} = 
[{ X}, { Y}]$ (torsion-freeness).\\

\noindent
Define the Christoffel symbol $\{\Gamma^k_{ij} \}$ by
$$
\nabla_{\frac{\partial}{\partial P_i}}
\frac{\partial}{\partial P_j}
=\sum^\ell_{k=1}
\Gamma^k_{ij} \frac{\partial}{\partial P_k}.
$$
Denote the $R$-module of K\"ahler differentials by $\Omega^{1}_{R}$.
Let 
\[
I^{*} : \Omega^{1}_{R} \times \Omega^{1}_{R} \longrightarrow R
\]
  be the symmetric $R$-bilinear form induced from $I^{*}$.  
Let $g^{ij}:=I^*(dP_i, dP_j)$ and 
$G:=(g^{ij})$, which is an $\ell \times \ell$-matrix with entries in $R$.
Note
\[ 
G = J({\bf P})^T A J({\bf P}).
\]
Define the contravariant Christoffel symbol $\{\Gamma^{ij}_k \}$ by
$$
\Gamma^{ij}_k:=-\sum^\ell_{s=1} g^{is} \Gamma^j_{sk}
$$
as in \cite[3.25]{dub1}. Define two $\ell \times \ell$-matrices
$$
\Gamma^*_k:=(\Gamma^{ij}_k), ~~~\Gamma_k:=(\Gamma^j_{ik}).
$$

\begin{lemma}
\label{2.2}
Let $1\leq k\leq \ell$. 

(1) $\Gamma^*_k=-G \Gamma_k,$

   (2) $\frac{\partial}{\partial P_k}[G]=\Gamma^*_k +(\Gamma^*_k)^T,$
   in particular,
   $D[G] = 
   \Gamma_{\ell}^{*}
   +
   (\Gamma_{\ell}^{*})^{T}, 
    $ 
    
   (3) $\Gamma^*_k=J({\bf P})^T A \frac{\partial}{\partial P_k}
[J({\bf P})],$
   
   (4) $\Gamma^*_\ell=B^{(1)}.$
\end{lemma}

\begin{proof} 
(1): By definition.   

(2): Apply \cite[3.26]{dub1}. 

(3): Let $S_k=(S^{ij}_k)$  be an $\ell\times\ell$-matrix defined by 
$$
S_k:=J({\bf P})^T A \frac{\partial}{\partial P_k}[J({\bf P})]. 
$$
It is enough to prove the (contravariant) compatibility (A) and the 
torsion-freeness (B) :

(A)  $\frac{\partial}{\partial P_k}[G]=S_k +S^T_k, $

(B)  $\sum_t g^{kt} S^{ij}_t=\sum_t g^{it} S^{kj}_t,$

\noindent
because of the uniqueness of the Levi-Civita connection.
We can verify (A) and (B) as follows : 

(A): 
\begin{eqnarray*} 
\pppk[G]
&=&
\pppk\left[
\JP^{T}A \JP 
\right]\\
&=&
\pppk\left[
\JP^{T}\right]
A \JP 
+
\JP^{T}
A \pppk\left[\JP\right] 
=
S_{k}^{T} + S_{k},  
\end{eqnarray*} 

(B):
\begin{eqnarray*} 
\sum_{t} g^{kt} S_{t}^{ij} 
&=&
\sum_{t, p, q, r, s}  
\frac{\partial P_{k}}{\partial X_{r} }
a^{rs}
\frac{\partial P_{t}}{\partial X_{s} }
\frac{\partial P_{i}}{\partial X_{p} }
a^{pq}
\frac{\partial}{\partial X_{t} }
\left(
\frac{\partial P_{j}}{\partial X_{q} }
\right)\\
&=&
\sum_{p, q, r, s}  
\frac{\partial P_{k}}{\partial X_{r} }
a^{rs}
\frac{\partial P_{i}}{\partial X_{p} }
a^{pq}
\sum_{t} 
\frac{\partial P_{t}}{\partial X_{s} }
\frac{\partial}{\partial X_{t} }
\left(
\frac{\partial P_{j}}{\partial X_{q} }
\right)\\
&=&
\sum_{p, q, r, s}  
\frac{\partial P_{k}}{\partial X_{r} }
a^{rs}
\frac{\partial P_{i}}{\partial X_{p} }
a^{pq}
\frac{\partial^{2}  P_{j}}{\partial X_{s} \partial X_{q} },
\end{eqnarray*} 
which is symmetric with respect to $i$ and $k$.
Thus
$
\sum_{t} g^{kt} S_{t}^{ij} 
=
\sum_{t} g^{it} S_{t}^{kj}. 
$  

(4) is easy, e. g., \cite[(2.17) (2.18)]{sot}.
\end{proof} 

\begin{remark}
If $P_{1}, \dots , P_{\ell}$ are chosen so that they satisfy the equality
\[
D[g^{ij}] = \delta_{i+j, \ell+1}, 
\]
it is known (e.g., \cite[pp. 275]{dub1}) that 
\[
\Gamma_{k}^{ij}/m_{j} 
=
\Gamma_{k}^{ji}/m_{i}. 
\]
By Lemma~\ref{2.2} (4) and (2), in this case, one has
\[
B^{(1)}_{ij} = \frac{m_{j}}{h} \delta_{i+j, \ell+1}.  
\]
By Lemma~\ref{2.1} (4),
\[
B^{(k)}_{ij} = \left\{(k-1) + \frac{m_{j}}{h}  \right\} 
\delta_{i+j, \ell+1}.
\]
\end{remark}

\begin{definition}
\label{formula}
Let $m \geq 0$.
Define $\xi^{(m)}_1, \dots, \xi^{(m)}_\ell \in \Der_S$ by:
$$
(\xi^{(m)}_1, \dots, \xi^{(m)}_\ell):= 
\begin{cases}
   (\frac{\partial}{\partial P_1}, \dots, \frac{\partial}{\partial P_\ell})
J({\bf P})^T AJ(D^k \left[{\bf X} \right])^{-1} \text{ if } m=2k, \\
   (\frac{\partial}{\partial P_1}, \dots, \frac{\partial}{\partial P_\ell})
J({\bf P})^T AJ(D^k \left[{\bf X} \right])^{-1} J({\bf P}) \text{ if } m=2k+1.
   \end{cases}
$$
\end{definition}

It is easy to see that

\begin{proposition}
\cite[Proposition 3.9]{ter}
\label{2.6}
For $k\geq 1$,
 \[
(\xi^{(2k+1)}_1, \dots, 
\xi^{(2k+1)}_\ell)=- (\xi^{(2k-1)}_1, \dots, \xi^{(2k-1)}_\ell)
(B^{(k)})^{-1} G.
\square
\]
\end{proposition}

The following result is the main theorem of \cite{ter}: 
\begin{theorem}
\cite[Theorem 1.1]{ter}
\label{2.5} 
Let $m \geq 0$.
 
(1)
The derivations $\xi^{(m)}_1, \dots, \xi^{(m)}_\ell$
form a basis for $\sD^{(m)}({\mathcal A}).$
In particular, 
$\sD^{(m)}({\mathcal A})$
is a free $S$- 
module of rank $\ell$,

(2) 
$$
\deg \xi^{(m)}_j=\begin{cases}
   kh ~~~(m=2k) \\
   kh+m_j ~~~(m=2k+1)
   \end{cases} 
   $$
 for $ 1 \leq j\leq \ell$. $\square$
\end{theorem}

\noindent
{\it Proof of Theorem \ref{2.4}.}  
It is easy to see that
each $\xi^{(2k+1)}_{j}  $ lies in $\Der_{R}$
by Definition \ref{formula}
and Proposition~\ref{2.6}. 
We have
\begin{eqnarray*} 
&~&
\left(
\nabla_{D} \xi^{(1)}_{1},
\dots,
\nabla_{D} \xi^{(1)}_{\ell}
\right)\\
&=&
\nabla_{D} \left(
\left(
\partial/\partial P_{1},
\dots,
\partial/\partial P_{\ell}
\right)
\JP^{T}A\JP 
\right)
=
\nabla_{D} \left(\left(
\partial/\partial P_{1},
\dots,
\partial/\partial P_{\ell}
\right)G
\right)\\
&=&
\left(
\nabla_{D} \left(
\partial/\partial P_{1}
\right),
\dots,
\nabla_{D} \left(
\partial/\partial P_{\ell}
\right)
\right)G
+
\left(
\partial/\partial P_{1},
\dots,
\partial/\partial P_{\ell}
\right) D[G]
\\
&=&
\left(
\partial/\partial P_{1},
\dots,
\partial/\partial P_{\ell}
\right)
\Gamma_{\ell}^{T} G 
+
\left(
\partial/\partial P_{1},
\dots,
\partial/\partial P_{\ell}
\right) (\Gamma_{\ell}^{*} + (\Gamma_{\ell}^{*})^{T}  )
\\
&=&
\left(
\partial/\partial P_{1},
\dots,
\partial/\partial P_{\ell}
\right) 
(- (\Gamma_{\ell}^{*})^{T}  + \Gamma_{\ell}^{*} + (\Gamma_{\ell}^{*})^{T}  )
\\
&=&
\left(
\partial/\partial P_{1},
\dots,
\partial/\partial P_{\ell}
\right) 
\Gamma_{\ell}^{*}
=
\left(
\partial/\partial P_{1},
\dots,
\partial/\partial P_{\ell}
\right) 
B^{(1)}
\end{eqnarray*} 
by Lemma \ref{2.2}. This shows (2) for $k=1$.

Next we will show (1) by induction on $k$.
By Propostion~\ref{2.6} and
Lemma~\ref{2.1}, one computes 
\begin{eqnarray*} 
&~&
\left(
\nabla_{D} \xi^{(2k+1)}_{1},
\dots,
\nabla_{D} \xi^{(2k+1)}_{\ell}
\right)
=
-\nabla_{D}\left( 
\left(
\xi^{(2k-1)}_{1},
\dots,
\xi^{(2k-1)}_{\ell}
\right)
(B^{(k)} )^{-1} G
\right)\\
&=&
-\left( 
\nabla_{D}\xi^{(2k-1)}_{1},
\dots,
\nabla_{D}\xi^{(2k-1)}_{\ell}
\right)
(B^{(k)} )^{-1}
 G
 \\
&~&~~~~
-
\left(
\xi^{(2k-1)}_{1},
\dots,
\xi^{(2k-1)}_{\ell}
\right)
(B^{(k)} )^{-1} D[G]\\
&=&
\left(
\xi^{(2k-3)}_{1},
\dots,
\xi^{(2k-3)}_{\ell}
\right)
(B^{(k-1)} )^{-1}
B^{(k)}  
(B^{(k)} )^{-1} 
G
 \\
&~&~~~~
-
\left(
\xi^{(2k-1)}_{1},
\dots,
\xi^{(2k-1)}_{\ell}
\right)
(B^{(k)} )^{-1}
(B^{(1)} + (B^{(1)})^{T} )\\
&=&
\left(
\xi^{(2k-3)}_{1},
\dots,
\xi^{(2k-3)}_{\ell}
\right)
(B^{(k-1)} )^{-1}
G
 \\
&~&~~~~
-
\left(
\xi^{(2k-1)}_{1},
\dots,
\xi^{(2k-1)}_{\ell}
\right)
(B^{(k-1)} )^{-1}
(B^{(k)} - B^{(k-1)})\\
&=&
\left(
\xi^{(2k-1)}_{1},
\dots,
\xi^{(2k-1)}_{\ell}
\right)
(-I_{\ell} - (B^{(k-1)})^{-1} B^{(k)} + I_{\ell})\\
&=&-
\left(
\xi^{(2k-1)}_{1},
\dots,
\xi^{(2k-1)}_{\ell}
\right)
(B^{(k-1)})^{-1} B^{(k)}.
\end{eqnarray*} 
Here, when $k=1$,
read $\xi_{j}^{(-1)} = \partial /\partial P_{j} 
~(j = 1, \dots , \ell)$ and $B^{(0)} = I_{\ell} $.   
Then this computation proves (1).
The assertion (2) for $k > 1$ easily 
follows from (1) and (2)
for $k=1$. $\square$  

\medskip

\begin{lemma}
\label{2.6B}
Let ${\mathcal T}_{k}:= 
\bigoplus_{j=1}^{\ell} T \xi^{(2k-1)}_{j}$   
for  $k\geq 1$. Then,
for $p\geq 1$, 
\[
\bigoplus_{k\geq p} {\mathcal T}_{k} 
=
\bigoplus_{j=1}^{\ell} R \xi^{(2p-1)}_{j}
=
\sD^{(2p-1)}(\A) \cap \Der_{R}.   
\]
\end{lemma}

\noindent
{\it Proof.}  
For the second equality, it is easy to see
\[
\bigoplus_{j=1}^{\ell} R \xi^{(2p-1)}_{j}
=
\left(\bigoplus_{j=1}^{\ell} S \xi^{(2p-1)}_{j}
\right)
 \cap \Der_{R} 
\]
by averaging over $W$.

Next we will prove the first equality.
Define a set 
${\mathcal S} := 
\{
\xi_{j}^{(2k-1)}
\mid
1\leq j\leq\ell, 1\leq k  
\}$.
First we will show that  
${\mathcal S} $ is linearly independent over 
$T$. 
Define matrices
$H_{k} ~~(k\geq 0)$ by
\[
H_{0} := I_{\ell} 
{\rm ~~~and~~~}  
H_{k} 
:= (-1)^{k} 
(B^{(1)})^{-1} G
(B^{(2)})^{-1} G
\dots
(B^{(k)})^{-1} G
~~~(k\geq 1).
\]
By applying Proposition \ref{2.6} inductively, we have 
\[
(\xi^{(2k+1)}_1, \dots, 
\xi^{(2k+1)}_\ell)=  (\xi^{(1)}_1, \dots, \xi^{(1)}_\ell)
H_{k}
\]
for $k\geq 0$. 
Since $D[B^{(k)} ] = 0$  and
$D^{2} [G] = D [B^{(1)} + \left(B^{(1)}\right)^{T}  ] = 0$ 
by Lemmas \ref{2.1} and \ref{2.2}, we obtain 
$D^{j} [H_{k} ] = 0$ for $j > k$. 
Suppose that ${\mathcal S} $ 
is linearly dependent over $T$. Then
there exist vectors $\mathbf{x}_{0}, \mathbf{x}_{1}, \dots , 
\mathbf{x}_{m}\in T^{\ell} $ 
satisfying
\[
{\mathbf x}_{m} \neq 0
{\rm ~~~and~~~} 
\sum_{k=0}^{m}  (\xi^{(2k+1)}_1, \dots, 
\xi^{(2k+1)}_\ell) {\mathbf x}_{k} = 0. 
\]
Thus we have
\begin{equation*} 
\sum_{k=0}^{m}  H_{k}  {\mathbf x}_{k} = 0. 
\end{equation*} 
Apply $D^{m} $ to both sides, and we get
$D^{m}[H_{m} ] {\mathbf x}_{m} = 0 $.
Note that $\det D[G]$ is a nonzero constant by K. Saito
\cite[5.1]{sai2}\cite[Corollary 4.1]{dub1}.  So
 \[
 D^{m}[H_{m} ] = (m!)
(B^{(1)})^{-1} D[G]
(B^{(2)})^{-1} D[G]
\dots
(B^{(k)})^{-1} D[G]
 \]
 is invertible.  Thus ${\mathbf x}_{m} = 0 $,
 which is a contradiction. This shows that 
 the set ${\mathcal S} $ is linearly independent over 
 $T$.  Therefore
 the sum
$
 \sum_{k\geq p} {\mathcal T}_{k}   
$
 is a direct sum.
  For $k\geq p$, 
  \[
  {\mathcal T}_{k}
  \subseteq
\sD^{(2k-1)}(\A) \cap \Der_{R} 
\subseteq 
\sD^{(2p-1)}(\A) \cap \Der_{R} 
=
\bigoplus_{j=1}^{\ell} R \xi^{(2p-1)} _{j}.
  \]
  Therefore we have
$$ 
\bigoplus_{k\geq p}{\mathcal T}_{k} \subseteq
\bigoplus_{j=1}^{\ell} R \xi^{(2p-1)} _{j}. 
$$
Compare the Poincar\'e series of both sides:
\begin{eqnarray*} 
{\rm Poin}(\bigoplus_{k\geq p}{\mathcal T}_{k}, t)
&=&
\left(\prod_{j=1}^{\ell-1} (1 - t^{m_{j}+1 } )
\right)^{-1} 
\sum_{k\geq p} 
\left(
t^{(k-1)h+m_{1} } +
\dots +
t^{(k-1)h+m_{\ell} } 
\right)\\
&=&
\left(
\prod_{j=1}^{\ell} (1 - t^{m_{j}+1 } )
\right)^{-1} 
\left(
t^{(p-1)h+m_{1} } +
\dots +
t^{(p-1)h+m_{\ell} } 
\right)\\
 &=&
{\rm Poin}
(\bigoplus_{j=1}^{\ell} R \xi^{(2p-1)} _{j}, t).
\,\,\,
\end{eqnarray*} 
This implies the first equality.
$\square$

\medskip

\noindent
{\it Proof of Theorem \ref{1.2}. }
It is enough to show $
{\mathcal G}_{k} 
=
{\mathcal T}_{k} 
 $ for each $k\geq 1$
 because of Lemma \ref{2.6B}.
 Since Theorem \ref{2.4} (2) asserts
$$
(\nabla^k_D \xi^{(2k-1)}_1, \dots, \nabla^k_D \xi^{(2k-1)}_\ell)
=(-1)^{k-1}(\frac{\partial}{\partial P_1}, \dots, \frac{\partial}
{\partial P_\ell}) B^{(k)},
$$
one has
\[
\nabla_{D}^{k} \xi^{(2k-1)}_{j} \in {\mathcal G}_{0}  \,\, 
(1\leq j\leq \ell, k\geq 1) 
\]
because $D[B^{(k)} ]=0$. 
Therefore we have  
$$ 
{\mathcal G}_{k} = \nabla_{D}^{-k} {\mathcal G}_{0} 
\supseteq
\bigoplus_{j=1}^{\ell} T \xi^{(2k-1)} _{j} = 
{\mathcal T}_{k}   
$$
for each $k\geq 1$. 
In \cite[6.3]{sai2} \cite[5.7]{sai4}, K. Saito showed 
\[
\bigoplus_{k\geq 1} {\mathcal G}_{k}
=
D_{R}(\Delta^{2} ) :=
\{
\theta \in \Der_{R} 
\mid
\theta(\Delta^{2} )\in \Delta^{2}R \}.  
\]
On the other hand, it is known \cite[Theorem 6.60]{ort} that
 \[
 D_{R}(\Delta^{2} )
 =
 \sD^{(1)}(\A) \cap \Der_{R}.
 \]
 Therefore, by Lemma \ref{2.6B}, we obtain
\[
\bigoplus_{k\geq 1}  {{\mathcal G}_{k}}
=
\bigoplus_{k\geq 1}  {{\mathcal T}_{k}}
\]
and thus
$
{\mathcal G}_{k} 
=
{\mathcal T}_{k} 
 $ for each $k\geq 1$.
 $\square$

\end{document}